\input amssym.def
\input amssym.tex
\magnification 1200

\def\newline{\break}
\def\item{\vskip .03 in}

\def\align{\eqalign}

\def\a{\alpha}

\def\un{\underline}

\def\pa{\partial}
\def\pr{\prime}
\def\noblackbox{\overfullrule=0pt}

\def\BZ{{ \Bbb Z}}
\def\BR{{ \Bbb R}}

\def\BC{{ \Bbb C}}

\def\lan{\langle}

\hoffset= 2pc
\hsize=30pc
\vsize=45pc

 \font\sectionfont=cmbx12 at 12pt



\def\sec#1{\vskip 1.5pc\noindent
{\hbox 
{{\sectionfont #1}}}\vskip 1pc}

\def\theo#1#2{\vskip 1pc
  {\bf   Theorem\ #1.} {\it #2}
\vskip 1pc}

\def\cor#1#2{\vskip 1pc{\bf 
\ Corollary\ #1.}
  {\it #2}\vskip 1pc}
\def\prop#1#2{\vskip 1pc{\bf 
 Proposition \ #1.} 
{\it #2}}

\def\proof#1{\vskip 1pc {\bf  
Proof.} {\rm #1} 
\vskip 1pc}

\def\qed{\vrule height4pt 
width3pt depth2pt}

\def\uBC*M{\u{\BC}^*~_M}

\def\uBC{\u{\BC}}

\def\text{\textstyle}

\def\text{\textstyle}

\def\Bi{Birkha\"user}

\def\2pii{2\pi\sqrt{-1}}

\def\display{\displaystyle}

\def\({\left(}
\def\){\right)}
\def\[{\left[}
\def\]{\right]}

\def\NSF92-95{\footnote{This research was supported in part by NSF grants
DMS-9203517 and DMS-9504522.}}
\def\NSF95{\footnote*{This research was supported in part by NSF grant
 DMS-9504522.}}
\def\NSF95{\footnote*{This research was supported 
in part by NSF grant
 DMS-9504522.}}
\centerline{{\bf The beta function of a knot}}
\vskip .5pc
\centerline{by Jean-Luc Brylinski\NSF95}
\vskip 1pc
\noblackbox
\sec{0. Introduction}

In this note we study a metric invariant of a knot $K$ embedded
in $\BR^3$. This metric invariant depends on a complex parameter
$s$. It is defined as a double integral using the parameterization
$x\mapsto \gamma(x)$ of $K$ by the arc-length parameter $x$:

$$B_K(s)=\int_{K\times K}~||\gamma(x)-\gamma(y)||^sdxdy\eqno(1)$$
where $dx$, $dy$ denote the arc-lengths on the two copies
of $K$. Clearly the function $(x,y)\mapsto ||\gamma(x)-\gamma(y)||^s$
is continuous for $Re(s)>0$, so that the integral converges
in that domain. One main result is

\theo{0.1}{For a  smooth knot $K$, the function
$s\mapsto B_K(s)$ extends analytically to a meromorphic function
on $\BC$, with only possible poles at $-1,-3,-5\cdots$.
The residue at $-2j-1$ is of the form
$\int_K~P_j(\kappa,\tau)dx$, where
$P_j$ is an explicitly computable polynomial
in the curvature $\kappa$, the torsion
$\tau$ and their derivatives.}

For $s=-1,-3,-5$ the polynomial $P_j$ is described
in \S 3.

The case $s=-2$ is $E(K)-4$, where $E(K)$ is the M\"obius energy
of $K$ studied in [{\bf F-H-W}] [{\bf O}] [{\bf K-K}] [{\bf K-S}].

The theorem is proved in section 3.

I would like to thank John Millson and Rob Kusner for
valuable comments on a first draft of this note, and in particular
for clarifying for me the relation with the M\"obius energy.

\sec{1. The case of the circular knot.}

We give the computation of $B_K(s)$ for the trivial
circular knot. In this case we have, using the rotational symmetry:
$$\align{B_K(s)&=2\pi \int_0^{2\pi}[2\sin ({x\over 2})]^sdx\cr
&=2^{s+2}\pi\int_0^{\pi}~sin(x)^sdx\cr
&=2^{s+3}\pi\int_0^{{\pi \over 2}}sin(x)^sdx}.$$
The integral is classically evaluated explicitly, via
the change of variable $sin(x)=\sqrt y$.
We have:
$$\int_0^{{\pi \over 2}}sin(x)^sdx={1\over 2}
\int_0^1y^{{s-1\over 2}}(1-y)^{-1\over 2}dy=B({s\over 2}
+{1\over 2},{1\over 2}),$$
using the Eulerian {\it beta function}
$$B(p,q)=\int_0^1~x^{p-1}(1-x)^{q-1}dx$$
which is equal to $\display {\Gamma(p)\Gamma(q)\over \Gamma(p+q)}$
[{\bf W-W}]. Therefore we find
$$\int_0^{{\pi \over 2}}|sin(x)|^sdx=
{\Gamma({1\over 2})\Gamma({(s+1)\over 2})\over
2\Gamma({s\over 2}+1)}$$

This is of course a meromorphic function of $s$,
and its  poles are located precisely at the poles
of $\Gamma({(s+1)\over 2})$, which occur exactly
when ${(s+1)\over 2}$ belongs to $\{ 0,-1,-2,\cdots\}$,
i.e., when $s\in \{ -1,-3,\cdots\}$.
This illustrates the theorem, and it also justifies calling
$B_K(s)$ the {\it beta function} of the knot $K$.

The functional equation for the gamma function immediately
gives the functional equation
$$B_K(s)=4{s-1\over s}B_K(s-2)$$
which can also be verified directly from the integral
expression for $B_K(s)$, using integration by parts.

\sec{2. The case of polygonal knots}

We consider a polygonal knot $K$, with $n$ sides $K_1,\cdots,K_n$
of lengths $a_1,\cdots,a_n$. It is convenient to think
of the set of edges as $\BZ/n\BZ$. We denote by $v_j$ the vertex
belonging to the sides $K_j$ and $K_{j+1}$ (so $v_n$ belongs
to $K_n$ and $K_1$). Let $\phi_j\in ( 0,\pi)$ be the angle
between these two edges. We study the function
$B_K(s)$ by decomposing it into ${n+1\over 2}$ double integrals
over $K_j\times K_l$, where $j<l$. If $K_j$ and $K_l$ do not meet, then
clearly the double integral $\int_{K_j\times K_l}~
||\gamma(x)-\gamma(y)||^sdxdy$
extends to a holomorphic function in $\BC$. So we have to
consider two types of double integrals:

$$A_j=\int_{K_j\times K_j}~||\gamma(x)-\gamma(y)||^sdxdy$$
and
$$B_j=\int_{K_j\times K_{j+1}}~||\gamma(x)-\gamma(y)||^sdxdy.$$

We have
$$A_j=\int_0^{a_j}\int_0^{a_j}~|x-y|^s
dxdy.$$
We evaluate this as follows: first we fix $x$ and integrate
over $y$, which yields the expression
$\display{1\over s+1}[x^{s+1}+(a_j-x)^{s+1}]$. So we have
$$A_j={2\over s+1}\int_0^{a_j}~x^{s+1}dx=
{2\over (s+1)(s+2)}a_j^{s+2},$$
which is meromorphic with simple poles at $s=-1$ and $s=-2$.
The residues are  $2a_j^{-1}$ and $-2$ respectively.

Now we study $B_j$. If $u$ resp. $v$ is the unit vector
pointing from $v_j$ to $v_{j-1}$ (resp. from $v_j$
to $v_{j+1}$), we have:
$$B_j=\int_0^{a_j}\int_0^{a_{j+1}}~||au-bv||^sdadb.$$
This is simplified somewhat by introducing the plane parallelogram
$P$ spanned by the vectors $a_ju$ and $-a_{j+1}v$.
We have in polar coordinates in the plane of $P$:
$$B_j=\int_P~sin(\phi_j)^{-1}r^srdrd\theta.$$
The factor $sin(\phi_j)^{-1}$ appears here as the jacobian
of the change of variables from $(a,b)$ to cartesian
coordinates $(x,y)$ in the plane
of $P$. Since we wish to determine
whether this function of $s$ has a meromorphic
extension and where the poles are located, we do not
 need to evaluate explicitly, in the sense that we 
may subtract the integral over any subregion which is clearly
holomorphic in $\BC$. Thus we may replace the parallogram
$P$ by its intersection with a disc of radius $\epsilon$
centered at $0$, to get the function
$$\align{C_j&=sin(\phi_j)^{-1}\int_0^{\epsilon}\int_0
^{\pi-\phi_j}~r^srdrd\theta\cr
&=sin(\phi_j)^{-1}(\pi-\phi_j){1\over s+2}[r^{s+2}]_0^{\epsilon}}$$
which for $Re(s)>-2$ is equal to
$${1\over s+2}\epsilon^{s+2}sin(\phi_j)^{-1}(\pi-\phi_j)$$
This makes it clear that $C_j$ is a meromorphic function
of $s$ with only a simple pole at $s=-2$, with residue equal
to $${\pi-\phi_j\over sin(\phi_j)}.$$
Now $B_K(s)-\sum A_j-2\sum C_j$ is an entire function.
Putting all this information together, we obtain

\prop{2.1}{(1) For a polygonal knot $K$, the beta function
$B_K(s)$ extends to meromorphic function of $s$, with at
 most simple poles at $s=-1$ and $s=-2$.
\vskip .04 in
(2) We have:
$$Res_{s=-1}~B_K(s)=l(K), Res_{s=-2}~B_K(s)=
-2n+2\sum_{j=1}^n~{\pi-\phi_j\over sin(\phi_j)}$$
where $l(K)$ is the length of $K$,
$n$ is the number of sides, and the $\phi_j$ are the angles
 between the two edges at the vertices $v_1,\cdots,v_n$.
The residues are always $>0$.}

The last statement follows from the fact that
$\pi-\phi_j>sin(\phi_j)$. We note that there is no harm
in introducing a fictitious vertex by splitting an edge
at some interior point. We note indeed that as $\phi$ tends
to $\pi$ the function ${\pi-\phi\over sin(\phi)}$ tends
to $1$. Then the residue at $s=-2$ is unchanged by this operation,
as it must be.

\sec{3. The beta function of a smooth knot.}

Clearly the problem with the double integral (1-1)
is that it will diverge near the diagonal $x=y$.
We write down the Taylor series for the parameterization
$x\mapsto \gamma(x)$ of $K$ by the arc-length parameter
$x$:
$$\gamma(y)-\gamma(x)=\sum_{j= 1}^r{1\over j!}
(y-x)^j\gamma^{(j)}(x)+(y-x)^{r+1}\a_r(x,y),\eqno(3-1)$$
where $\a_r$ is a smooth vector-valued function
of $(x,y)$. For the purpose of studying the analytic properties
of $B_K(s)$, we may replace it by
the function
$$D(s)=\int_{|x-y|\leq \epsilon}~||\gamma(x)-\gamma(y)||^sdxdy
\eqno(3-2)$$
for some small $\epsilon$, which will be determined later.
For the square of the distance function we have an 
expression of the type:
$$||\gamma(x)-\gamma(y)||^2=|y-x|^2(1+\sum_{i=1}^{r-1}~f_i(x)(y-x)^i
+(y-x)^{r}Q(x,y)),$$
where each $f_i$ is a smooth function of $x$, and $Q(x,y)$
is also smooth. When we take the $s\over 2$-th power, we use
the power series expansion
$(1+a)^{s\over 2}=\sum_{j\geq 0}~{{s\over 2}\choose j}a^j$ which
has radius of convergence $1$, and gives a holomorphic
function of $(a,s)$ in the domain $s\in\BC,|a|<1$. We now assume that
$\epsilon$ is picked small enough so that
$|\sum_{i=1}^{r-1}~f_i(x)(y-x)^i
+(y-x)^{r}Q(x,y)|<1$ whenever $|y-x|<\epsilon$.
Then in that region
we obtain an expansion:
$$||\gamma(x)-\gamma(y)||^s=|y-x|^s(1+\sum_{i=1}^{r-1}h_i(s,x)(y-x)^i
+(y-x)^{r}R(s,x,y)),$$
where $h_i(s,x)$ is a smooth function of $(s,x)$ which
is polynomial as a function of $s$, and $R(s,x,y)$
is a smooth function of $(s,x,y)$ which is an  entire
function of $s$.
We then have
$$D(s)=\int_0^l\int_{-\epsilon}^{\epsilon}~
\bigl[|y-x|^s+\sum_{i=1}^{r-1}h_i(s,x)|y-x|^s(y-x)^i
+R(s,x,y)|y-x|^s(y-x)^{r}dxdy\eqno(3-3)$$
In the region $Re(s)>-r$ the last term of the integral
is holomorphic, so the behaviour of $D(s)$ can be inferred
form the other terms. The integral over $y$ is easily evaluated,
and it is clear  that for $i$ odd it is equal to $0$, whereas
for $i=2$ even we have 
$$\int_{-\epsilon}^{\epsilon}~|y-x|^{s+2j}={2
\over s+2j+1}\epsilon^{s+2j+1}.$$
This is meromorphic in $\BC$, with a simple
pole at $s=-2j-1$ with residue equal to
$2$. It then follows that $D(s)$ is meromorphic
in the region $Res(s)>-r$ with possible poles
at $-1,-3,\cdots,-r+1$, with residue
at $-2j-1$ equal to
$2\int_0^l~h_j(-2j-1,x)dx.$
Since this is true for all $r\geq 1$, Theorem 0.1
is proved. The fact that $h_j(s,x)$ has a differential expression
in terms of $\kappa$ and $\tau$ follows since any scalar product
$\gamma^{(p)}(x)\cdot \gamma^{(q)}(x)$ of derivatives of $\gamma$
has such an expression.

We now compute the residues at $s=-1$ and $s=-3$.
We easily obtain the following expression for the
square of the distance function:

$$\align{||\gamma(x)-\gamma(y)||^2&=(y-x)^2
(1+(y-x)^2({1\over 8}
||\gamma^{\pr\pr}(x)||^2
+{1\over 3}\gamma{\pr}(x)\cdot \gamma{\pr\pr\pr}(x)+O(y-x)^4)\cr
&=(y-x)^2
(1-{5\over 24}(y-x)^2||\gamma^{\pr\pr}(x)||^2+O(y-x)^4)}$$
using the fact that $\gamma^{\pr}\cdot \gamma^{\pr\pr}=0$
and the corollary that
$\gamma^{\pr}\cdot \gamma^{\pr\pr\pr}+\gamma^{\pr\pr}\cdot
\gamma^{\pr\pr}=0$.
We then have
$$||\gamma(x)-\gamma(y)||^s=|y-x|^s(1-{5s\over 48}||\gamma^{\pr\pr}(x)||^2
(y-x)^2+O(y-x)^4).$$
We then derive the formula of the first two residues:
$$Res_{s=-1}B_K(s)=2l=2l(K)~,~Res_{s=-3}B_K(s)={5\over 8}
\int_K\kappa^2dx.$$
The residue at $-5$ can be similarly calculated. It is equal
to $\int_KQdx$, where
$$Q={(\kappa^{\pr})^2\over 8}
-{\kappa^2\tau^2\over 144}+{859\over 16\times 144}\kappa^4.$$

We now compare the value $B_K(-2)$ of the beta function
with the M\"obius energy $E(K)$  studied 
by  O'Hara [{\bf O}] and by  Freedman, He and Wang [{\bf F-H-W}]. This energy is defined as

$$E(K)=\int_{K\times K}~[~||\gamma(x)-\gamma(y)||^{-2}-d(x,y)^{-2}]dx dy,$$
where $d(x,y)=\int_x^yds$ is the distance between $x$ and $y$
along the knot.

This functional is the value at $s=-2$ of the following
functional

$$D_K(s)=\int_{K\times K}~[||\gamma(x)-\gamma(y)||^{s}-d(x,y)^{s}] dx dy.$$

This function of $s$ is clearly holomorphic in the region
$Re(s)>-3$ from the asymptotic expansion of $||\gamma(x)-\gamma(y)||^{s}$
previously described. The function
$$f(s)=\int_{K\times K}d(x,y)^{s}dx dy$$
is easily evaluated to be equal to $\display{2^{-s}l^{s+2}\over s+1}$, which
is holomorphic near $s=-2$ with value at $s=-2$ equal to $4$.

Hence we find

\prop{3.1}{The functional $B_K(-2)$ is equal to $E(K)-4$.}

In particular, the funtional $B_K(-2)$ coincides with the version
of the M\"obius energy given in [{\bf K-K}], which is normalized
to vanish for the circle knot. J. Sullivan pointed out to me
that this vanishing, in the case of $B_K(-2)$, follows immediately
from the functional equation given in section 1.

\sec{4. Poisson brackets}

In [{\bf Br}] we introduced a symplectic structure
over the Fr\'echet manifold $Y$ of all smooth knots in
$\BR^3$. This induces a Poisson bracket over a class
of smooth functionals $f$, called
supersmooth. This class of functionals
$f$ is such that the differential $df$ at a knot $K$ can be viewed
as a section of the restriction to $K$ of the cotangent bundle
of $\BR^3$, which we may view as a vector field to $\BR^3$
defined along $K$. This vector field is denoted by $df$. We then have
$$\{ f,g\}_K=-l(K)^-2\int_Kdet({d\gamma\over dx},df(x)
,dg(x))dx$$ [{\bf Br, (3-11), p. 138}]. We apply this to two functionals
$K\mapsto B_K(s)$ and $K\mapsto B_K(u)$. These two functionals
will be denoted by $B_s$ and $B_u$. We first assume that the real
parts of $s$ and $u$ are large enough. This will allow us to
compute as many derivatives of the functional $K\mapsto B_s(K)$ as needed.
We see that the derivative
$dB_s$ at a knot $K$ is given by the vector-valued function on $K$:
$$[dB_s]_K(x)=-2s\int_K~||\gamma(x)-\gamma(y)||^{s-1}
({\gamma(x)-\gamma(y)})dy.$$
Therefore the Poisson bracket is defined and given by the 
triple integral

$$\align{\{ B_s,B_u\}_K=4sul(K)^{-2}\int_{K\times K\times K}~
&||\gamma(x)-\gamma(y)||^{s-1}||\gamma(x)-\gamma(z)||^{u-1}\cr
&det({d\gamma\over dx},\gamma(x)-\gamma(y),\gamma(x)-\gamma(z))
dxdydz.}$$
Now, while keeping $Re(u)$ large, we want to move $s$
to a vicinity of $s=-1$. We should then expand the integrand
in powers of $y-x$. Since we have
$$det({d\gamma\over dx},\gamma(x)-\gamma(y),\gamma(x)-\gamma(z))
=(x-y)^2det({d\gamma\over dx},\gamma^{\pr\pr}(x),
\gamma(x)-\gamma(z))+O(x-y)^3$$
in a neighborhood of the partial diagonal $y=x$,
it is clear that the integrand
is of the form $|y-x|^{s+1}$ times a continuous function.
Therefore the Poisson bracket $\{ B_s,B_u\}$ has no
singularity at $s=-1$. This implies

\prop{4.1}{The length functional $l$ on the
space of knots $Y$ Poisson commutes
with all the functionals $u\mapsto B_K(u)$.}

\cor{4.2}{For any $j\geq 1$, the length
functional on $Y$ Poisson commutes with
the functional $K\mapsto Res_{s=-2j-1}~B_K(s)$.}

The fact that $l$ Poisson commutes with the functional
$\int_K~\kappa^2dx$ is well-known (see [{\bf Br}]).

\sec{5. A Bernstein-Sato type equation.}

For a polynomial $f(z_1,\cdots,z_n)$, J. Bernstein [{\bf Be}]
proved that there exists a partial differential
equation
$$P(\un z,{\pa\over \pa \un z},s)f^s=B(s)f^{s-1},\eqno(5-1)$$
where $P$ is a differential operator with polynomial coefficients,
which depends polynomially in the variable $s$.
Here $B(s)$ is a monic polynomial; when it is chosen of minimal
degree, it is called the Bernstein-Sato polynomial of $f$
and the equation (5-1) is called a Bernstein-Sato equation.
In some cases the equation was known for quite some time:
for instance, if $f(z_1,\cdots,z_n)=z_1^2+\cdots+z_n^2$
and if we set $P(\un z,{\pa\over \pa \un z})
={1\over 4}\sum_{i=1}^n~{\pa^2\over \pa z_i^2}$, then
we have
$$Pf^s=s(s-1+{n\over 2})f^{s-1}.$$
This equation can be found in the treatise of
Gelfand and Shilov [{\bf G-S}]. 

We will establish an analog of this equation for the functional
$B_K(s)$ on the space fo knots $Y$. First we write down an analog
of the Laplace operator. Given a knot $K$ and $2$ points of $K$
with arc-length parameters $x,y$, we introduce
the unit vector
$$w(x,y)={\gamma(x)-\gamma(y)/||\gamma(x)-\gamma(y)||}.$$
This is of course a smooth vector-valued function on $K\times K$.

Now let $f$ be a supersmooth functional on $Y$.
Then its second order differential $D^2f$, evaluated at
$K$, is an $\BR^3\otimes \BR^3$-valued function
on $K\times K$.
We set
$$P(f)(K)=\int_{K\times K}~\lan D^2f(x,y),w(x,y)\otimes w(x,y)
dxdy\eqno(5-2)$$

We then have

\prop{5.1}{For $Re(s)>>0$, we have:
$$PB_s=s(s-1)B_{s-2}\eqno(5-3)$$}

\proof{This follows immediately from the fact that
the second derivative of $||\gamma(x)-\gamma(y)||^s$
evaluated at $(x,y)$ on the element
$w\otimes w$ is equal
to $s(s-1)||\gamma(x)-\gamma(y)||^{s-2}$. \qed}
\vskip .03 in
\vskip 1.5 pc
\centerline{REFERENCES}

\vskip .7 pc
[{\bf Be}] J. Bernstein, \it The analytic continuation
of generalized functions with respect to a parameter,
\rm Funct. Anal. Appl. \bf  6\rm (1972), 272-285

[{\bf Br}] J-L. Brylinski, \it Loop Spaces,
Characteristic Classes and Geometric Quantization,
\rm Progress in Math. vol. \bf 107 \rm (1993), \Bi

[{\bf F-H-W}] M.H. Freedman, Z.X. He and Z. Wang, 
\it M\"obius energy of knots and unknots,
\rm Annals of Math \bf 139 \rm (1994), 1-50

[{\bf G-S}] I. M. Gel'fand and G. E. Shilov,
\it Generalized Functions, \rm Academic Press 

[{\bf K-K}] D. Kim and R. Kusner, 
\it Torus knots extermizing the M\"obius energy,
\rm Experiment. Math. \bf 2\rm (1993), 1-9

[{\bf K-S}]  R. B. Kusner and J. M. Sullivan,
    \it {M\"obius} Energies for Knots and Links, Surfaces
            and Submanifolds, \rm
     in Geometric Topology,
   ed. W.  H. Kazez,
    Amer. Math. Soc./Int'l Press (1997), 570-604

[{\bf O}] J. O'Hara, \it Energy of a Knot,
\rm Topology \bf 30 \rm (1991), 241-247

[{\bf W-W}] E. T. Whittaker and G. N. Watson, 
\it A course in Modern Analysis, \rm
Cambridge Univ. Press

\vskip .12 in
Penn State University

Department of Mathematics

305 McAllister

University Park, PA. 16801

USA

\vskip .03 in
e-mail:jlb@math.psu.edu

\bye